\documentclass{elsart1}

\usepackage{amssymb}
\usepackage{mathrsfs}
\usepackage{amsmath}
\newcommand{\R}{{\Bbb R}}

\newcommand{\N}{{\Bbb N}}

\begin{document}
\begin{frontmatter}
\title{On the geometry of wave solutions of a
delayed reaction-diffusion equation}
\author[a]{Elena Trofimchuk,}
\author[b]{Pedro Alvarado}
\author[b]{and Sergei Trofimchuk}
\address[a]{Department of Diff. Equations, National Technical University, Kyiv,  Ukraine
 \\ {\rm E-mail: trofimch@imath.kiev.ua} }
\address[b]{Instituto de Matem\'atica y Fisica, Universidad de Talca, Casilla 747,
Talca, Chile \\ {\rm E-mails: pedro.alvarado.r@gmail.com,
trofimch@inst-mat.utalca.cl}}

\medskip
\begin{center}
\it To Yuriy Alekseevich Mitropolskiy
\end{center}
\medskip
\end{frontmatter}

\section{Introduction and main results} \label{intro}

\vspace{-5mm}

This paper concerns positive wave solutions of the non-local delayed
reaction-diffusion equation
\begin{equation}\label{17} \hspace{-7mm}
u_t(t,x) = u_{xx}(t,x)  - u(t,x) + \int_{\R}K(x-y)g(u(t-h,y))dy, \ u
\geq 0,
\end{equation}
which is widely used in applications, e.g. see
\cite{Brit,fhw,gouk,gouS,gouss,LRW,MSo,tz,wlr}. In part, our
research was inspired by several problems  raised in
\cite{gouss,LiWu,ma1,swz}. We suppose that equation (\ref{17}) has
exactly two equilibria $u_1 \equiv 0, \ u_2 \equiv \kappa >0$ and

\vspace{-3mm}

\begin{eqnarray}\label{kkk}\hspace{-7mm}
K \geq 0, \ \ \int_{\R}K(s)e^{\lambda s}ds \ {\rm is \ finite \ for
\ all} \ \lambda \in \R, \ \ {\rm and}\ \int_{\R}K(s)ds = 1.
\end{eqnarray}
Note that the usual restriction  $K(s) = K(-s), \ s \in \R,$ is
not required here. In a biological context, $u$ is the size of an
adult population, so we will consider  only non-negative solutions
of equation (\ref{17}). The nonlinear $g$ is called {\it the birth
function}, it is assumed to satisfy the following hypothesis
\begin{description}
\item[{\rm \bf(H)}] $g \in C(\R_+, \R_+)$ has only one local
extremum at $s=s_M$ (maximum) and   $g(0)=0$, $g(s)
> 0$ if $s > 0$. Next, $0$ and $\kappa >0$ are the only two solutions  of $g(s)=s$, and $g$ is differentiable at $s=0$, with
$g'(0)>1$.
\end{description}
For example, this is the case in the Nicholson's blowflies model
\cite{fhw,gouss,LRW,ma1,MSo} where $g(s)= pse^{-s}$. See also
Subsection 1.6 below.

Let us  fix some terminology. Following \cite{GK}, we call bounded
positive classical  solutions $u(x,t) = \phi(x +ct)$ satisfying
$\phi(-\infty) = 0$ semi-wavefronts. We say that the
semi-wavefront $u(x,t) = \phi(x +ct)$  is a wavefront [is a
pulse], if the profile function $\phi$ satisfies $\phi(+\infty) =
\kappa$ [respectively, satisfies $\phi(+\infty) = 0$].  Wavefronts
are the most studied subclass of semi-wavefronts. Asymptotically
periodic semi-wavefronts represent another subclass, see
\cite{TT}. Some of our results are proved for semi-wavefronts, and
some of them, for  wavefronts. For example, the setting of
semi-wavefronts is more convenient to work with the problem of the
minimal speed of propagation, cf. \cite[Section 3]{tz}.

In Subsections 1.1-1.5 below, we present our main results. Their
proofs and some additional comments can be found in Sections 2-9.

{\bf 1.1 Two critical speeds and non-existence of pulse waves.} In
this subsection, we consider more general equation
\begin{equation}\label{eG}
u_t = u_{xx}  - q u + F(u, \mathcal{K}_1u, \dots ,
\mathcal{K}_mu),
\end{equation}
where $F:\R_+^{m+1} \to \R_+$ is a continuous function and
$$(\mathcal{K}_ju)(t,x): = \int_{\R}K_j(x-y)f_j(u(t-h,y))dy.$$
This equation includes (\ref{17}) as a particular case (as well as
equations (\ref{17al}), (\ref{WW}) considered below).  We assume
that each kernel $K_j$ satisfies condition (\ref{kkk}) and the
continuous non-negative functions $F, f_j$ are differentiable at
the origin. Set
\begin{equation}\label{pqK}
p:= \sum_{j=1}^m  F_{s_j}(0)f_j'(0),  \ K(s): = p^{-1}\sum_{j=1}^m
F_{s_j}(0)f_j'(0)K_j(s).
\end{equation}
We assume that $p
> q >0$ and $F(0) = f(0)= F_{s_0}(0) =0;$ and $F_{s_j}(0),
f_j'(0) \geq 0$  for all $j =1,\dots,m$. As a consequence,  $K$
satisfies (\ref{kkk}). Next, consider
\begin{equation}\label{har}
 \psi(z,\epsilon) =\epsilon
z^2-z-q+p\exp(-zh)\int_ {\R} K(s)\exp(-\sqrt{\epsilon} zs)ds,
\end{equation}
and let $\epsilon_i =\epsilon_i(h,p,q)> 0, i =0,1,$ be as in Lemma
\ref{L23} of Appendix. Set $c_*:= 1/\sqrt{\epsilon_0}$ and $c_\#:=
1/\sqrt{\epsilon_1}$. By Lemma \ref{L23},  $c_* \geq c_\#$, and
$c_* = c_\#$ if and only if $c_* = c_\#=0$.  As we show in
Appendix, $c_\# =0$ if $\int_{\R}sK(s)ds \geq 0$ and $c_*
> 0$ if $\int_{\R}sK(s)ds \leq 0$. Moreover,  $c_\# >0$ if the equation
$$z^2 -q + p\int_{\R}\exp(-zs)K(s)ds=0$$ has
negative roots. The main result of this subsection  is the
following
\begin{thm} \label{41a22}  Let $u(t,x) = \phi(x+ct), \ c > 0,$ be a positive
bounded solution of equation (\ref{eG}).  If $c < c_*$ then $
\liminf\limits_{s \to -\infty} \phi(s)
> 0$ and therefore $\phi(x+ct)$ is
not a semi-wavefront. Next, if $c > c_\#$, then $\phi(x+ct)$ is
persistent: $ \liminf\limits_{s \to +\infty} \phi(s) > 0.$ In
consequence, equation (\ref{eG}) does not have non-stationary
pulses.
\end{thm}
Observe that even when $g$ is monotone on $[0,\kappa]$, it was not
known whether every semi-wavefront  to equation (\ref{17}) is
separated from zero as $x+ct \to +\infty$. The persistence of
semi-wavefronts was established in \cite{ETST} for a local version
of model (\ref{17}). The proof in \cite{ETST} is based on the
local estimations technique which does not apply to equation
(\ref{17}). To overcome this obstacle, we will use a Laplace
transform approach developed in \cite{FA,FA2} and successfully
applied in \cite[Proposition 4]{gous}, \cite[Theorem 4.1]{pituk},
\cite[Theorem 5.4]{ETST}.

We emphasize that $c_*$ can be different from the minimal speed of
propagation of semi-wavefronts even for the simpler case of
equation (\ref{17}), cf. \cite{GK}. However, as it was shown in
\cite{tz,ETST,WLR}, $c_*$ coincides with the minimal speed for Eq.
(\ref{17}) if $g$ satisfies {\rm \bf(H)} together with the
additional condition
\begin{equation}\label{SL}
g(s) \leq g'(0)s \quad {\rm for\ all \ } \quad s \geq 0.
\end{equation}
\begin{rem}
A lower bound for the admissible speeds of semi-wavefronts to the
reaction-diffusion functional equation
$$u_t(t,x) = \Delta
u(t,x) + g(u_t), \ u(t,x) \geq 0,\ x \in \R^n,
$$
was already calculated in the pioneering work  of Schaaf, see
Theorem 2.7 (i) and Lemma 2.5 in \cite{sch}. Recent work
\cite{ETST} complements  Schaaf's investigation in two aspects:
(i) analyzing the case of non-hyperbolic trivial equilibrium and
(ii) taking into account the problem of {\it small solutions}.

Note that very few theoretical studies are devoted to the minimal
speed problem for the  {\it non-local} equation (\ref{eG}). To the
best of our knowledge, the first accurate proof of the
non-existence of semi-wavefronts was provided by Thieme and Zhao
in \cite[Theorem 4.2 and Remark 4.1]{tz} for the equation
\begin{equation}\label{17al}
u_t = u_{xx}  - f(u) + \int_{\R}K_{\alpha}(x-y)g(u(t-h,y))dy, \
K_{\alpha}(x) = \frac{e^{-x^2/(4\alpha)}}{\sqrt{4\pi\alpha}}.
\end{equation}
In order to prove this result, Thieme and Zhao have extended an
integral-equations approach \cite{die,tie} to scalar non-local and
delayed reaction-diffusion equations. Their proof makes use of the
special form of the kernel $K$ which is the fundamental solution of
the heat equation.

Besides the above mentioned work \cite{tz}, a non-existence result
was proved for the equation
\begin{equation}\label{WW}
u_t(t,x) = u_{xx}(t,x)  + g(u(t,x), \int_{\R}K(x-y)u(t-h,y)dy),
\end{equation}
in the recent work \cite{WLR} by Wang, Li and Ruan. Their method
required $C^2$-smoothness of $g$ and the fulfillment of several
convexity conditions.

Our approach is different from those in \cite{tz} and \cite{WLR}
and it allows us to impose minimal restrictions on the right-hand
side of equation (\ref{eG}). In any case, the problem of
non-existence of semi-wavefronts to equations (\ref{eG}),
(\ref{17al}), (\ref{WW}) is non-trivial, and the corresponding
proofs are not easy. In fact, some papers provide only a heuristic
explanation for why non-local models similar to (\ref{eG}) do not
have positive wavefronts propagating at velocity $c$ which is less
than some critical speed $c_*$;  see, for instance,
\cite{ZAMP,gouss,LRW,swz,wlr}. In the mentioned works, $c_*$ is
defined as the unique positive number for which some associated
characteristic function $\psi_{c_*}$ (similar to (\ref{har})) has
a positive multiple root while $\psi_c$ does not have any positive
root for all $c < c_*$, cf. Lemma \ref{L23}. However, it seems
that this argument is incomplete. Indeed, some linear autonomous
functional differential equations of mixed type may have a
nonoscillatory solution in spite of the nonexistence of real roots
of its characteristic equation. See remarkable examples proposed
by Krisztin  in \cite{kr}.
\end{rem}
{\bf 1.2. Uniform persistence of waves} The second aspect of the
problem we address is the uniform persistence of positive waves
$u(t,x)= \phi(x +ct), \ c > c_\#,$ to Eq. (\ref{17}). This
property means that $\liminf_{s \to +\infty} \phi(s)\geq \zeta$,
where $\zeta>0$ depends only on $g$. The uniform persistence of
positive bounded waves will be proved by assuming  condition
(\ref{kkk}) and the following hypothesis
\begin{description}
\item[{\rm \bf(B)}] $g \in C(\R_+,\R_+)$ satisfies $g(s)
>0$ when $s
> 0$ and, for some $0 <\zeta_1 < \zeta_2$,

1. $g([\zeta_1,\zeta_2])\subseteq [\zeta_1,\zeta_2]$ and
$g([0,\zeta_1])\subseteq [0,\zeta_2]$;

2. $\min_{s \in [\zeta_1,\zeta_2]}g(s) = g(\zeta_1)$;

3.  (i) $g(s) > s$ for $s \in (0, \zeta_1]$ and (ii) there exists
$p=g'(0) \in (1, +\infty)$;

4. In $\R_+$, the equation $g(s)=s$ has exactly two solutions $0$
and $\kappa$.
\end{description}
Remark that conditions in {\rm \bf(B)} are weaker than {\rm \bf(H)}.
Indeed, set $\zeta_2 = g(s_M)$  if $g$ satisfies {\rm \bf(H)}. It is
easy to see that the map $g: [0,\zeta_2] \to [0,\zeta_2]$ is well
defined. We can also consider the restrictions $g: [\zeta_1,\zeta_2]
\to [\zeta_1,\zeta_2]$ for every positive $\zeta_1 \leq \min \{g
^2(s_M),\ s_M\}$. Clearly, there exists $\zeta_1$ satisfying {\rm
\bf(B.2)}
 \vspace{-2mm}
\begin{thm} \label{mainA} Assume {\rm \bf(B)} and let $u= \phi(x +ct), \ c > c_\#,$ be a
positive bounded solution to equation (\ref{17}) . Then $ \zeta_1
\leq \liminf_{s \to +\infty} \phi(s) \leq \sup_{s \in \R} \phi(s)
\leq \sup g([0,\sup_{s \in \R} \phi(s)]). $
\end{thm}
{\bf 1.3 Existence of semi-wavefronts} Set $\tilde c_* :=
(\epsilon_0(h,\sup_{s >0}g(s)/s,1))^{-1/2}$. It is clear that
$\tilde c_* \geq c_* $ and  $\tilde c_*  = c_* $ if conditions
(\ref{SL}) and {\rm \bf(B-3ii)} hold. Our third result establishes
the existence of semi-wavefronts for all $c \geq \tilde c_*$:
 \vspace{-2mm}
\begin{thm} \label{mainAA} Assume {\rm \bf(B)} except the condition {\rm \bf(B-3ii)}.
Then equation (\ref{17}) has a positive semi-wavefront $u(t,x)=
\phi(x +ct)$  for every positive $c \geq \tilde c_*$.
\end{thm}
 \vspace{-2mm}
The proof of Theorem \ref{mainAA} relies on the Ma-Wu-Zou method
proposed in \cite{wz} and further  developed in \cite{ma,ma1,swz}.
It uses the positivity and monotonicity properties of the integral
operator
 \vspace{-3mm}

\begin{eqnarray} \label{n} \hspace{-9mm}
  (A\phi)(t)&=& \frac{1}{\epsilon'}\left\{\int_{-\infty}^te^{\lambda (t-s)}(G\phi)(s-h))ds +
\int_t^{+\infty}e^{\mu (t-s)}(G\phi)(s-h))ds \right\},\nonumber  \\
 \hspace{-9mm}  & & (G\phi)(s)  =
 \int_{\R}K(w)g(\phi(s-\sqrt{\epsilon}w))dw,\ \epsilon':=\epsilon (\mu -
\lambda)
\end{eqnarray}

 \vspace{-3mm}

where $\lambda < 0 < \mu$ solve $\epsilon z^2 -z -1 =0$ and
$\epsilon^{-1/2} = c > 0$ is the wave velocity. As it can be
easily observed, the profiles $\phi \in C(\R,\R_+)$ of travelling
waves are completely determined by the integral equation
$A\phi=\phi$ and the Ma-Wu-Zou method  consists in the use of an
appropriate fixed point theorem to $A:\mathfrak{K }\to
\mathfrak{K}$, where  $\mathfrak{K} = \{x: 0\leq \phi^-(t) \leq
\phi(t) \leq \phi^+(t)\}$ is subset of an adequate Banach space
$(C(\R,\R), |\cdot|)$. Now, $\mathfrak{K}$ should be 'nice' enough
to assure the compactness (or monotonicity) of $A$. This
requirement is not easy to satisfy. Thus only relatively narrow
subclasses of $g$ (e.g. sufficiently smooth at the positive
equilibrium and monotone or quasi-monotone in the sense of
\cite{wz}) were considered within this approach. Our contribution
to the above method is the very simple form of the bounds $\phi^{\pm}$
for $\mathfrak{K}$. For instance, due to the information provided
by Theorem \ref{mainA}, we may take $\phi^-(t) =0$ for all $t \geq
0$.  Here, this finding allows to weaken the smoothness conditions
imposed on $g(s)$ at $s =0$. In particular, for equation
(\ref{17}), Theorem \ref{mainAA} improves Theorem 1.1 in
\cite{ma1}. Indeed, the method employed in \cite{ma1} needs
essentially that $K(s) = K(-s)$ and $\limsup_{u\to +0}(g'(0) -
g(u)/u)u^{-\nu}$ is finite for some $\nu \in (0,1]$.

{\bf 1.4 Delay-depending conditions of the existence of
wavefronts} As it happened in the case of semi-wavefronts, the
existence of the wavefronts depends not only on the derivatives
$g'(0),g'(\kappa)$ but also on the values of the entire function
$g$. Here, we prove their existence by analyzing some
one-dimensional dynamical systems associated to $g$. The property
of the negative Schwarzian $(Sg)(s)=g'''(s)/g'(s)-(3/2)
\left(g''(s)/g'(s)\right)^2$ is instrumental in simplifying the
analysis of these systems in some cases, see  Proposition \ref{SW}
and further comments in Appendix. The next result presents
delay-depending conditions of the existence of the wavefronts, it
follows from more general Theorem \ref{mainex}.
\begin{thm} \label{mainex2} Assume {\rm \bf(H)} and let $\zeta_1, \zeta_2$ be
as in {\rm \bf(B)}. Suppose further that $(Sg)(s)<0$,  $s \in
[\zeta_1,\zeta_2]\setminus\{s_M\}$ and $g^2(\zeta_2) \geq \kappa$.
If, for some $\epsilon >0$,
$$(1-\min\{e^{-h}, \int_{-h/\sqrt{\epsilon}}^{0}K(u)du\})g'(\kappa) \geq -1,$$  then
Eq. (\ref{17}) has a wavefront $u(t,x) =\phi(x +ct)$ for every $c
\geq \max\{\tilde c_*,1/\sqrt{\epsilon}\}$. Moreover, for these
values of $c$, each semi-wavefront is in fact a wavefront .
\end{thm}
{\bf 1.5 Non-monotonicity of wavefronts} The problem of
non-monotonicity of wavefronts to equation (\ref{17}) was widely
discussed in the literature. The state of the art is surveyed in
\cite[Section 4.3]{gouss}. As far as we know, the paper
\cite{ZAMP} by Ashwin {\it et al.} contains the first heuristic
explanation  of this phenomenon.  Recent works \cite{FT,TT} have
provided rigorous analysis  of non-monotonicity in the local case.
Here, we follow the approach of \cite{FT} to indicate conditions
inducing the loss of monotonicity of wavefronts in the simpler
case  when the kernel $K$ has compact support.

\begin{thm} \label{main} Suppose $g$ is continuous, $g(\kappa)=\kappa,\
g'(\kappa)<0,$  supp $K \subseteq [-\eta,\eta]$, and the equation
\begin{equation}\label{cchh1} (z/c)^2-z-1+g'(\kappa)\exp(-zh)\int_
{-\eta}^{\eta} K(s)\exp(-zs/c)ds =0
\end{equation} does not have any root in $(-\infty,0)$ for some fixed $c=\bar{c}$.
If $\phi(+\infty) = \kappa$ for a non-constant solution
$\phi(x+\bar{c}t)$ of equation (\ref{17}), then  $\phi(s)$
oscillates about $\kappa$.
\end{thm}
{\bf 1.6 An example}\label{7} We apply our results to the
reaction-diffusion-advection equation
\begin{equation}\label{RADE}
\frac{\partial u}{\partial t} = \frac{\partial }{\partial
x}(D_m\frac{\partial u}{\partial x} + Bu) - u +
\int_{\R}K_1(x+Bh-y)g(u(t-h,y))dy,
\end{equation}
where $g(s) = pse^{-s}, \ K_\alpha(s) =
(4\pi\alpha)^{-1/2}e^{-s^2/(4\alpha)}$. This equation was studied
numerically in \cite{LiWu} for various values of parameters $p, B,
h, D_m$. Plugging the traveling wave Ansatz $u(x,t) = \phi(x+ct)$
into (\ref{RADE}), we obtain that
$$
D_m\phi''(t)-(c-B)\phi'(t)-\phi(t) + \int_{\R}K(s)g(\phi(t-(c-B)h
-s))ds = 0.
$$
Next, setting $\gamma = 1/D_m, \  \phi(s) = z(s/(c-B)), \ \epsilon
= D_m/(c-B)^2$, we find that
\begin{equation*}
\epsilon z''(t)- z'(t)-z(t) +
\int_{\R}K_{\gamma}(s)g(z(t-h-\sqrt{\epsilon}s))ds = 0.
\end{equation*}
For this equation, $g$ satisfies condition (\ref{SL}) and the
function $\psi(z,\epsilon)$ from (\ref{har}) can be found
explicitly:
$\psi(z,\epsilon)  = \epsilon z^2 -z-1+p\exp(\epsilon\gamma z^2-z
h).$
An easy calculation shows that $\kappa = \ln p, $ $ \zeta_2 = p/e, $
$  g'(\kappa) = \ln(e/p)$.  As in Section 4.1 of \cite{LiWu}, we
select $D_m = 5, \ h = 1, \ p = 9$. Then we find that $g^2(\zeta_2)=
3.299 > \kappa= 2.197$. Analyzing $\psi(z,\epsilon)$, we obtain that
$\epsilon_0= 0.3725\dots$. In consequence, equation (\ref{RADE}) has
semi-wavefronts if and only if $c-B \geq \sqrt{5/0.3725\dots} =
3.66\dots$ This can explain (see also Remark \ref{re}) the emergence
of unsteady multihump waves in the numerical experiments realized in
\cite{LiWu}: indeed, the value $c-B = 3$ taken in \cite{LiWu} is
less than the minimal speed of semi-wavefronts. Finally, an
application of Theorem \ref{mainex2} shows that equation
(\ref{RADE}) with $D_m = 5, \ h = 1, \ p = 9$ has positive {\it
wavefronts} if and only if $c - B \geq 3.66\dots$

\vspace{-5mm}

\section{Proof of Theorem \ref{41a22} for $c < c_*$}\label{fi}

\vspace{-5mm}

 Let $u(t,x)= \phi(ct+ x)$ be a positive bounded
solution of (\ref{eG}) and suppose that $\epsilon := c^{-2} >
\epsilon_0(h,p,q)$. Set $\varphi(s)=\phi(cs)$. Then $\xi(t)=
\varphi(-t)$ satisfies
\begin{equation}\label{twe2a}
\epsilon \xi''(t) + \xi'(t)- q\xi(t)+(\mathcal{F}\xi)(t) =0, \quad
t \in \R,
\end{equation}
where   $(\mathcal{F}\xi)(t) = F(({\mathscr{\bf I}}\xi)(t))$ with
$({\bf I}\xi)(t)\in \R_+^{m+1}$ denoting
$$(\xi(t),
\int_{\R}K_1(s)f_1(\xi(t+\sqrt{\epsilon}s+h))ds,\dots,
\int_{\R}K_m(s)f_m(\xi(t+\sqrt{\epsilon}s+h))ds).$$ Since $\xi(t)$
is a bounded solution of equation (\ref{twe2a}), it must satisfy
\begin{equation}\label{iie2}
\xi(t) = \frac{1}{\epsilon (\tilde{\mu }-
\tilde{\lambda})}\left\{\int_{-\infty}^te^{\tilde{\lambda}
(t-s)}(\mathcal{F}(\xi))(s)ds + \int_t^{+\infty}e^{\tilde{\mu}
(t-s)}(\mathcal{F}(\xi))(s)ds \right\},
\end{equation}
where $\tilde{\lambda} < 0 < \tilde{\mu}$ are roots of $\epsilon
z^2+z-q=0$.

The following inequality is crucial in the coming discussion.
\begin{lem} \label{la} If $\xi:\R \to (0,+\infty)$ is
a bounded solution of equation (\ref{twe2a}),  then
\begin{equation}\label{inest}
\xi(t) \geq e^{\tilde{\lambda} (t-s)}\xi(s), \ t \geq s.
\end{equation}
\end{lem}

\vspace{-7mm}

\begin{pf} Since $(\mathcal{F}\xi)(t)$ is non-negative, after differentiating (\ref{iie2}),
we obtain
$$
\xi'(t)-\tilde{\lambda} \xi(t) =
\frac{1}{\epsilon}\int_t^{+\infty}e^{\mu (t-s)}(\mathcal{F}\xi)(s)ds
\geq 0.
$$
Therefore $(\xi(t)e^{-\tilde{\lambda} t})' \geq 0$, which implies
(\ref{inest}). \hfill {$\blacksquare$}
\end{pf}

Arguing by contradiction, we suppose that $\liminf\limits_{t \to
+\infty} \xi(t) =0$ for some  $\epsilon> \epsilon_0$.

\underline{\textbf{Case I: }} $\limsup\limits_{t \to +\infty}
\xi(t) = \lim\limits_{t \to +\infty} \xi(t) =0$.\\ In virtue of
(\ref{inest}) and Lemma \ref{41aa} from Appendix, we can find a
real number $D
>1$ and a sequence $t_n \to +\infty$ such that $\xi(t_n) = \max_{s
\geq t_n}\xi(s)$ and
$$
\max_{s \in [t_n-  8\sqrt{\epsilon},t_n]}\xi(s) \leq D \xi(t_n).
$$
It is easy to see that, for every fixed $n$, $\xi'(t)$ is  either
negative on $[t_n-  4\sqrt{\epsilon},t_n]$  or there is $t'_n \in
[t_n- 4\sqrt{\epsilon},t_n]$ such that $\xi'(t'_n) =0$, and $\xi(t)
\leq \xi(t'_n)$ for all $t \geq t_n'$ (thus $\xi(t) \leq D\xi(t'_n)$
if $t \in [t'_n- 4\sqrt{\epsilon},t'_n]$). Now, if $\xi'(t)$ is
negative then
$$|\xi'(t''_n)| =(\xi(t_n-4\sqrt{\epsilon}) -
\xi(t_n))/(4\sqrt{\epsilon}) \leq (D-1)\xi(t_n)/(4\sqrt{\epsilon})
:= D_1\xi(t_n)$$ for some $t_n'' \in [t_n- 4\sqrt{\epsilon},t_n]$.
Since $\xi(t''_n) \geq \xi(t_n)$, we obtain that
$$
|\xi'(t''_n)|  \leq  D_1\xi(t_n) \leq D_1 \xi(t''_n), \ {\rm and} \
\xi(t) \leq D\xi(t''_n) \ {\rm for \ all } \ t \in[t''_n-
4\sqrt{\epsilon},t_n''].
$$
Hence, by the above reasoning,  we may assume that $D$ and
$\{t_n\}$ are such that
$$
|\xi'(t_n)| \leq D \xi(t_n), \ \max_{s \in
[t_n-4\sqrt{\epsilon},t_n]}\xi(s) \leq D\xi(t_n)\ {\rm  and}\ \xi(t)
\leq \xi(t_n), \ t \geq t_n.
$$
Next, since  continuous $F$ is differentiable at $0$ and $F(0)=0$,
we obtain that
$$
F(s_0,s_1,\dots,s_m) = \sum_{j=0}^mA_j(s_0,s_1,\dots,s_m)s_j, \ s_j
\geq 0,
$$
where $A_j$ are continuous and $A_j(0) = F_{s_j}(0), \ j=1,\dots,m$.
In consequence, $y_n(t) = \xi(t+t_n)/\xi(t_n), t \in \R,$ should
satisfy the equation
\begin{equation}\label{twerin2}  \hspace{-8mm}
\epsilon y''(t) + y'(t)-a_{0,n}(t)y(t)+
\int_{\R}\mathscr{K}_{n}(t,s,\epsilon)y(t+\sqrt{\epsilon}s+h)ds=0,
\end{equation}
where $\mathscr{K}_{n}(t,s,\epsilon):=
\sum_{j=1}^mK_j(s)a_{j,n}(t,t+\sqrt{\epsilon}s+h)$ and
$$
a_{0,n}(t):= q-A_0(({\bf I}\xi)(t+t_n)), \  a_{j,n}(t,u):=
A_j(({\bf I}\xi)(t+t_n))\frac{f_j(\xi(u+t_n))}{\xi(u+t_n)}.
$$
From (\ref{inest}), it is clear that $e^{\tilde{\lambda}t} \leq
y_n(t) \leq 1$ for all $t \geq 0$ and $y_n(t)\leq
e^{\tilde{\lambda}t}, \ t \leq 0$. In particular, $y_n(0)=1$. Note
also that there is $C_{\xi}
>0$ such that $|a_{j,n}(t,u)| \leq C_{\xi} $ for all $j=0,\dots,m;
n \in \N; t,u \in \R$. Moreover, $\lim_{n \to \infty} a_{0,n}(t)=
q,$ $\lim_{n \to \infty} a_{j,n}(t,u) = F_{s_j}(0)f_j'(0)$
pointwise. Set $$\mathscr{G}_n(t):=
\int_{\R}\mathscr{K}_{n}(t,s,\epsilon)y_n(t+\sqrt{\epsilon}s+h)ds.$$
We claim that for arbitrary fixed $\sigma, \tau
>0$ there exists $c_{\sigma, \tau} >0$ such that $|\mathscr{G}_n(t)| \leq c_{\sigma, \tau}$
for all $t\in [-\sigma, \tau]$ and for all $n \in \N$.  Indeed,
$$
|\mathscr{G}_n(t)| \leq
\int_{\R}|\mathscr{K}_{n}(t,s,\epsilon)|y_n(t+\sqrt{\epsilon}s+h)ds
\leq
C_{\xi}\int_{\R}\sum_{j=1}^mK_j(s)y_n(t+\sqrt{\epsilon}s+h)ds\leq
$$
$$
\leq
C_{\xi}\int_{-\frac{t+h}{\sqrt{\epsilon}}}^{+\infty}\sum_{j=1}^mK_j(s)ds
+C_{\xi}\int^{-\frac{t+h}{\sqrt{\epsilon}}}_{-\infty}\sum_{j=1}^mK_j(s)e^{\tilde{\lambda}(t+\sqrt{\epsilon}s+h)}ds\leq
$$
$$
\leq
C_{\xi}(m+e^{\tilde{\lambda}(h-\sigma)}\int_{\R}\sum_{j=1}^mK_j(s)e^{\tilde{\lambda}\sqrt{\epsilon}s}ds)=:c_{\sigma,
\tau}, \ t\in [-\sigma, \tau].
$$
Now, since $z_n(t) = y'_n(t)$ solves the initial value problem
$z_n(0) = \xi'(t_n)/\xi(t_n) \in [-D,D]$ for
$$
\epsilon z'(t) + z(t)-a_{0,n}(t)y_n(t)+ \mathscr{G}_n(t) =0,
$$
we deduce the existence of  $k_{\sigma, \tau}>0$ such that, for
all $t\in [-\sigma, \tau]$ and $n \in \N$,
\begin{equation}\label{nv}
|y'_n(t)| = |e^{-t/\epsilon}z_n(0) +
\frac{1}{\epsilon}\int_{0}^te^{(s-t)/\epsilon}(a_{0,n}(s)y_n(s)-
\mathscr{G}_n(s))ds| \leq
\end{equation}
$$
\leq De^{\sigma/\epsilon}+
\frac{1}{\epsilon}|\int^{t}_0e^{(s-t)/\epsilon}(C_{\xi}\max\{1,
e^{-\tilde{\lambda}s}\}+ c_{\sigma, \tau})ds| \leq  k_{\sigma,
\tau}.
$$
Therefore, we may apply the Ascoli-Arzel\'a compactness criterion
together with a diagonal argument on each of the intervals
$[-i,i]$ to find a subsequence $\{y_{n_j}(t)\}$  converging, in
the compact-open topology, to a non-negative function $y_*:\R \to
\R_+$. It is evident that $y_*(0) = 1$ and $e^{\tilde{\lambda}t}
\leq y_*(t) \leq 1$ for all $t \geq 0$ and $y_*(t)\leq
e^{\tilde{\lambda}t}, \ t \leq 0$. By the Lebesgue's dominated
convergence theorem, we have that,  for every $t \in \R$,
$$ \mathscr{G}_n(t)
\to \mathscr{G}_*(t):= p\int_{\R}K(s)y_*(t+\sqrt{\epsilon}s+h)ds,
$$
where $K(s), p$ are defined in (\ref{pqK}). In consequence,
integrating (\ref{nv})  without $|\cdot|$ between $0$ and $t$ and
then taking the limit as $n_j \to \infty$ in the obtained
expression, we establish that $y_*(t)$ satisfies
\begin{equation}\label{eqw}
\epsilon y''(t) + y'(t)-qy(t)+
p\int_{\R}K(s)y(t+\sqrt{\epsilon}s+h)ds=0.
\end{equation}
Then Lemma \ref{41lin} implies that $ y_*(t) =  w(t) +
O(\exp(2\tilde{\lambda} t)), \ t \to +\infty, $ where $w$ is a non
empty finite sum of  eigensolutions of (\ref{eqw}) associated to
the eigenvalues $\nu_j \in F= \{2\tilde{\lambda} < \Re{\nu}_j \leq
0\}$. Observe now that $\nu$ is an eigenvalue of (\ref{eqw}) if
and only if $-\nu $ is a root of (\ref{har}). In this way, $F$
does not contain any real eigenvalue for $\epsilon > \epsilon_0$
(by Lemma \ref{L23}), and therefore $y_*(t)$ should be oscillating
on $\mathbb{R}_+$, a contradiction.

\underline{\textbf{Case II: }} $\liminf\limits_{t \to +\infty}
\xi(t) =0 \ {\rm and} \ S= \limsup\limits_{t \to +\infty} \xi(t)
> 0.$\\
In case II, for every fixed $j > S^{-1}$ there exists a sequence
of intervals $[p'_i,q'_i]$, $ \lim p'_i  = +\infty$ such that
$\xi(p'_i)= 1/j, \ \lim \xi(q'_i) =0, \ \xi'(q'_i)= 0$ and $\xi(t)
\leq 1/j,$ $t \in [p'_i,q'_i]$. Note that $\limsup(q'_i-p'_i) =
+\infty$ since otherwise we get a contradiction: the sequence
$\xi(t+p'_i)$ of solutions to equation (\ref{iie2}) contains a
subsequence converging to a non-negative bounded solution
$\xi_1(t)$ such that $\xi_1(0) = 1/j, \ \xi_1(\sigma) = 0$ for
some finite $\sigma
> 0$. In consequence, $w_i(t) = \xi(t+p'_i),$ $ \ t \in \R$ has a
subsequence converging to some bounded non-negative solution
$w_*(t)$ of (\ref{iie2}) satisfying $0 < w_*(t) \leq 1/j$ for all
$t \geq 0$. Since the case $w_*(+\infty) =0$ is impossible due to
the first part of the proof, we conclude that $0 < S^* =
\limsup\limits_{t \to +\infty} w_*(t) \leq 1/j$. Let $r_i \to
+\infty$ be such that $w_*(r_i) \to S^*$, then $w_*(t+r_i)$ has a
subsequence converging to a positive solution $\zeta_j: \R \to [0,
1/j]$ of (\ref{iie2}) such that $\max_{t \in \R}\zeta_j(t) =
\zeta_j(0)= S^* \leq 1/j$. Next, arguing as in case I after
formula (\ref{twerin2}), we can use sequence
$\{y_j(t):=\zeta_j(t)/\zeta_j(0)\}$ to  obtain a bounded positive
solution $y_*: \R \to (0,1)$ of linear equation (\ref{eqw}). For
the same reason as given in  Lemma \ref{la}, {\it bounded} $y_*$
decays at most exponentially. Now, invoking Lemma \ref{41lin} and
the oscillation argument as in case I, we get a contradiction.

\vspace{-3mm}
\section{Proof of Theorem \ref{41a22} for $c > c_\#$}
The case $c > c_\#$ is similar to case considered in Section
\ref{fi}. Below we give some details. Let $u(t,x)= \phi(ct+ x)$ be
a positive bounded solution of (\ref{eG}) and suppose that
$\epsilon := c^{-2} < \epsilon_1(h,p,q)$. Set
$\varphi(s)=\phi(cs)$ and $\tilde{K}_j(s)=
{K}_j(-s-2h/\sqrt{\epsilon})$. Then each $\tilde{K}_j(s)$
satisfies (\ref{kkk}) and $\varphi(t)$ verifies
$$
\epsilon \varphi''(t) - \varphi'(t)-
q\varphi(t)+(\mathcal{H}\xi)(t) =0, \quad t \in \R, $$ where
$(\mathcal{H}\varphi)(t) = F(({\mathscr{\bf J}}\varphi)(t))$ with
$({\bf J}\varphi)(t)\in \R_+^{m+1}$ denoting
$$(\varphi(t),
\int_{\R}\tilde{K}_1(s)f_1(\varphi(t+\sqrt{\epsilon}s+h))ds,\dots,
\int_{\R}\tilde{K}_m(s)f_m(\varphi(t+\sqrt{\epsilon}s+h))ds).$$
Since $(\mathcal{H}\xi)(t)$ is non-negative, the same argument as
used to prove Lemma \ref{la} shows that $\varphi(t) \geq
e^{\lambda (t-s)}\varphi(s), \ t \geq s,$ where
${\lambda}=-\tilde{\mu} < 0 < {\mu}=-\tilde{\lambda}$ are the
roots of $\epsilon z^2-z-q=0$. All this  allows to repeat the
proof given in Section \ref{fi}, with a few obvious changes, to
establish the persistence of $\varphi(t)$. For example, the
paragraph below (\ref{eqw}) should be modified in the following
way:

 "\dots we establish that $y_*(t)$ satisfies
\begin{equation}\label{eqw0}
\epsilon y''(t) - y'(t)-qy(t)+
p\int_{\R}K(s)y(t+\sqrt{\epsilon}s+h)ds=0.
\end{equation}
Then Lemma \ref{41lin} implies that $ y_*(t) =  w(t) +
O(\exp(2\lambda t)), \ t \to +\infty, $ where $w$ is a {\it non
empty} finite sum of  eigensolutions of (\ref{eqw0}) associated to
the eigenvalues $\lambda_j \in F= \{2\lambda <  \Re\lambda_j \leq
0\}$. Now, since the set $F$ does not contain any real eigenvalue
for $\epsilon \in (0,\epsilon_1)$ (see Lemma \ref{L23}), we
conclude that $y_*(t)$ should be oscillating on $\mathbb{R}_+$, a
contradiction".
\vspace{-5mm}
\section{Proof of Theorem \ref{mainA}}
Let $\varphi(t)$ satisfy $0 < \varphi(t)\leq M_0, \ t \in \R,$ and
\begin{equation}\label{twe}
\epsilon \varphi''(t) - \varphi'(t)-\varphi(t)+
\int_{\R}K(s)g(\varphi(t-\sqrt{\epsilon}s-h))ds=0, \quad t \in \R.
\end{equation}

\vspace{-5mm} Being bounded, $\varphi$ must verify the integral
equation
\begin{equation}\label{iie}
\varphi(t) = \frac{1}{\epsilon'}\left\{\int_{-\infty}^te^{\lambda
(t-s)}(G\varphi)(s-h)ds + \int_t^{+\infty}e^{\mu
(t-s)}(G\varphi)(s-h)ds \right\},
\end{equation}
where  $\epsilon',\mu, \lambda$ and $G\varphi$ are as in
(\ref{n}). From (\ref{iie}), we obtain that $ |\varphi'(t)| \leq
\max_{s \in [0,M_0]}g(s)/\epsilon'.$ This implies the
pre-compactness of the one-parametric family $\mathcal{F} =
\{\varphi(t+s), s \in \R\}$ in the compact open topology of
$C(\R,\R)$. It is an easy exercise to prove (by using (\ref{iie}))
that the closure of $\mathcal{F}$ consists from the positive
bounded solutions of (\ref{twe}). Next, for $\varphi$ as above,
set
$$0 \leq m = \inf\limits_{t \in \R}\varphi(t) \leq \sup\limits_{t \in \R}\varphi(t) = M < +\infty.$$
\begin{lem} \label{vg} $[m,M] \subseteq g([m,M])$.
\end{lem}
\vspace{-5mm}
\begin{pf} Indeed, if  $M=\varphi(s')=\max_{s \in \R}\varphi(s)$, a straightforward
estimation  of the right hand side of (\ref{iie}) at $t = s'$
generates ${M} \leq \max_{{m} \leq s \leq {M}} g(s)$. As long as
the maximum $M$ is not reached,  using the pre-compactness of
$\mathcal{F}$, we can find a solution $z(t)$ of $(\ref{twe})$ such
that $z(0) = \max_{s \in\R}z(s) = M$ and $ \inf_{s \in\R}z(s) \geq
m$. Therefore, by the above argument, ${M} \leq \max_{{m} \leq s
\leq {M}} g(s)$. The inequality ${m} \geq \min_{{m} \leq s \leq
{M}} g(s)$ can be proved in a similar way. Thus we can conclude
that $[m,M] \subseteq g([m,M])$. \hfill {$\blacksquare$}
\end{pf}
\vspace{-5mm} Note that Lemma \ref{vg} implies that
$\sup\varphi(t) \leq \sup g([0,\sup\varphi(t)])$.

Analogously, we have
\begin{lem} \label{dud} Let $\varphi$ satisfy (\ref{twe}) and be such that
$0 \leq m' = \liminf\limits_{t \to +\infty}\varphi(t) \leq
\limsup\limits_{t \to +\infty}\varphi(t) = M' < +\infty.$ Then
$[m',M'] \subseteq g([m',M'])$.
\end{lem}
\begin{thm} \label{43} Assume $\mathbf{(B)}$ and consider a positive bounded solution $\varphi\not\equiv 0$
of equation (\ref{twe}) for some fixed $\epsilon\in
(0,\epsilon_1)$. If ${m} = \inf_{s \in \R}\varphi(s) < \zeta_1$
then, in fact, $\epsilon\in (0,\epsilon_0)$ and $ \lim_{t \to
-\infty} \varphi(t) =0.$
\end{thm}
\vspace{-5mm}
\begin{pf} Set ${M}=\sup_{s \in \R}\varphi(s)$, then Lemma \ref{vg}
guarantees that $[{m},{M}] \subseteq g([{m},{M}])$. The
assumptions $\mathbf{(B)}$  and ${m} < \zeta_1$ make impossible
the inequality ${m}
> 0$. In consequence, ${m} = 0$ and, due to Theorem \ref{41a22},
either $\varphi(-\infty) = 0$ or $$0 = \liminf\limits_{t \to
-\infty} \varphi(t) < \limsup\limits_{t \to -\infty} \varphi(t) =
S. $$ However, as we will show it in the continuation, the  second
case cannot occur. Indeed, otherwise for every positive $\delta_1
<\min\{\zeta_1,S\}$, it would be possible to indicate two
sequences of real numbers $p_n < q_n$ converging to $-\infty$ such
that $\varphi(p_n) = \max_{[p_n,q_n]}\varphi(u) =\delta_1$, and
$\varphi(q_n)< \varphi(s)< \varphi(p_n)$ for all $s \in (p_n,q_n)$
with $ \lim \varphi(q_n) = 0$. We notice that necessarily
$\lim(q_n-p_n) = +\infty$, since in the opposite case an
application of the "compactness argument" leads to the following
contradiction: the sequence of solutions $\varphi(t+p_n)$ contains
a subsequence converging to a solution $\psi \in C(\R,\R)$ of
equation (\ref{iie}) verifying $\psi(0) = \delta_1$ and $\psi(t_0)
= 0$, for some finite $t_0 > 0$. Hence, $\lim(q_n-p_n) = +\infty$
and the limit solution $\psi$ is positive and such that $\psi(0) =
\delta_1 = \max_{s \geq 0} \psi(s)$. Moreover, by Theorem
\ref{41a22}, we have that $\delta_0:=\liminf\limits_{t \to
+\infty} \psi(t) >0$. In consequence, using again the "compactness
argument", we can construct a solution $\tilde{\psi}(t)$ of
equation (\ref{iie}) such that $\delta_0 \leq \tilde{\psi}(t) \leq
\delta_1 < \zeta_1,$ for all $t \in \R$. But, in view of
hypotheses $\mathbf{(B)}$, this contradicts to Lemma \ref{vg}.
\hfill {$\blacksquare$}
\end{pf}

\vspace{-4mm}

Now we are ready to prove that $ \liminf\limits_{t \to +\infty}
\varphi(t) \geq \zeta_1. $ Indeed, otherwise,  by the "compactness
argument", we can construct a bounded solution
$\tilde{\varphi}(t)$ such that $0 <\liminf\limits_{t \to +\infty}
\varphi(t) \leq \inf_{s \in \R}\tilde{\varphi}(s) < \zeta_1,$
contradicting to Theorem \ref{43}.

\vspace{-2mm}

\section{An application of Ma-Wu-Zou reduction}

\vspace{-2mm}

\noindent Throughout this section, $\chi_{\R_-}(t)$ stands for the
indicator of $\R_-$. Following the notations of Lemma \ref{L23} in
Appendix, for given $ \epsilon \in (0, \epsilon_0)$ we will denote
by $\lambda_1 = \lambda_1(\epsilon)< \lambda_2 =
\lambda_2(\epsilon)$ the positive roots of $\psi(z,\epsilon) =0$.
Also we will require
\begin{description}
\item[{\rm \bf(L)}] $g: (0, +\infty) \to (0, +\infty)$ is bounded
and locally linear in some right $\delta$-neighborhood of the
origin: $g(s) = ps, \ s \in [0,\delta)$, with $p > 1$.
Furthermore, $g(s) \leq ps$ for all $s \geq 0$.
\end{description}
Assuming this, for every $ \epsilon \in (0, \epsilon_0)$, we will
prove the existence of semi-wavefronts of  equation (\ref{twe}).
As it was shown by Ma, Wu and Zou \cite{ma,ma1,swz,wz}, solving
(\ref{twe}) can be successfully reduced to the determination of
fixed points of the integral operator $A$ from (\ref{n}) which is
considered in some closed, bounded, convex and $A$- invariant
subset $\mathfrak{K}$ of an appropriate Banach space $(X,
\|\cdot\|)$. In this section, the choice of $\mathfrak{K} \subset
X$ is restricted by the following natural conditions: (i) constant
functions cannot be elements of $X$; (ii) the convergence
$\varphi_n \to \varphi$ in $\mathfrak{K}$ is equivalent to the
uniform convergence $\varphi_n \Rightarrow \varphi_0$ on compact
subsets of $\R$. With this in mind, for some $\rho \in
(\lambda_1,\mu)$ and $\delta$ as in {\rm \bf(L)}, we set
\vspace{-4mm}
\begin{eqnarray*}
  X &=& \{\varphi \in C(\R,\R): \|\varphi\| =
\sup_{s \leq 0} e^{-\lambda_1s/2 }|\varphi(s)|+ \sup_{s \geq 0}
e^{-\rho s}|\varphi(s)|<
\infty\}; \\
 \mathfrak{ K }&=& \{\varphi \in X: \phi^-(t)= \delta(e^{\lambda_1 t} -
e^{\lambda_2 t})\chi_{\R_-}(t) \leq \varphi(t) \leq \delta
e^{\lambda_1t} = \phi^+(t), \ t \in \R\}.
\end{eqnarray*}
A formal linearization of $A$ along the trivial steady state is
given by
$$
(L\varphi)(t) = \frac{p}{\epsilon'}\left\{\int_{-\infty}^te^{\lambda
(t-s)}(\mathrm{Q}\varphi)(s-h)ds + \int_t^{+\infty}e^{\mu
(t-s)}(\mathrm{Q}\varphi)(s-h)ds \right\},
$$
where
$$
(\mathrm{Q}\varphi)(s) =
\int_{\R}K(w)\varphi(s-\sqrt{\epsilon}w)dw.
$$
\begin{lem} We have
$L\phi^{+} = \phi^{+}$. Next, \ $(L\psi)(t)> \psi(t), \ t \in \R,$
where $$\psi(t) := (e^{\lambda_1 t} - e^{\nu t})\chi_{\R_-}(t) \in
\mathfrak{K}$$ is considered with $\nu \in (\lambda_1,
\lambda_2]$.
\end{lem} \vspace{-2mm}
\begin{pf} It suffices to prove that $(L\psi)(t)> \psi(t)$ for $t \leq 0$. But we have
$$
(L\psi)(t)> \frac{p}{\epsilon'}\left\{\int_{-\infty}^te^{\lambda
(t-s)}(\mathrm{Q}(e^{\lambda_1(\cdot)} - e^{\nu(\cdot)}))(s-h)ds +
\right.
$$
$$\left. \int_t^{+\infty}e^{\mu (t-s)}(\mathrm{Q}(e^{\lambda_1(\cdot)} -
e^{\nu(\cdot)}))(s-h)ds \right\} \geq \psi(t).
$$
\hfill {$\blacksquare$}
\end{pf} \begin{lem} Let assumption {\rm \bf(L)} hold and
$\epsilon \in (0, \epsilon_0)$. Then  $A(\mathfrak{K}) \subseteq
\mathfrak{K}$.
\end{lem}

\vspace{-5mm}

\begin{pf} We have $A\varphi \leq L\varphi \leq L\phi^+ = \phi^+$ for every $\varphi \leq
\phi^+$. Now, if for some $u = s - \sqrt{\epsilon}w$ we have $0
<\phi^-(u) \leq \varphi(u)$, then $u <0$, so that $\varphi(u) \leq
\delta e^{\lambda_1u} < \delta$ implying $g(\varphi(u)) =
p\varphi(u) \geq p\phi^-(u)$. If $\phi^-(u_1) =0$ then again
$g(\varphi(u_1)) \geq p\phi^-(u_1)=0$. Therefore $(G\varphi)(t) \geq
p(\mathrm{Q}\phi^-)(t), \ t \in \R$ so that $A\varphi \geq L\phi^-
> \phi^-$ for every $\varphi \in \mathfrak{K}$. \hfill
{$\blacksquare$}
\end{pf}

\vspace{-5mm}

\begin{lem}$\mathfrak{K}$ is a closed, bounded, convex subset of $X$
and $A:\mathfrak{K }\to \mathfrak{K}$ is completely continuous.
\end{lem}

\vspace{-4mm}

\begin{pf} Note that the convergence of a sequence in $\mathfrak{K}$
amounts to the uniform convergence on compact subsets of $\R$. Since
$g$ is bounded, we have $|(A\varphi)'(t)|\leq \max_{s \geq
0}g(s)/\epsilon'$ for every $\varphi \in \mathfrak{K}$. The lemma
follows now from the Ascoli-Arzel${\rm \grave{a}}$ theorem combined
with the Lebesgue's dominated convergence theorem. \hfill
{$\blacksquare$}
\end{pf} \begin{thm} \label{34} Assume $({\bf L})$ and let $\epsilon \in (0, \epsilon_0)$.
Then the integral equation (\ref{iie}) has a positive bounded
solution in $\mathfrak{K}$.
\end{thm}

\vspace{-4mm}

\begin{pf} Due to the above lemmas, we can apply the Schauder's fixed point
theorem to $A:\mathfrak{K} \to \mathfrak{K}$. \hfill
{$\blacksquare$}
\end{pf}

\section{Proof of Theorem \ref{mainAA} }

\vspace{-4mm}

\underline{\bf Case I: $c > \tilde{c}_*$.}  First, we assume that
$\max_{s \geq 0} g(s) = \max_{s \in [\zeta_1,\zeta_2]} g(s) \leq
\zeta_2$. Set $k= \sup_{s
>0}g(s)/s$ (so that  $ks \geq g(s)$ for all $s \geq 0$) and
consider the  sequence
$$
\gamma_n(s)=\left\{%
\begin{array}{ll}
    ks, & \hbox{ for} \ s \in [0,1/(nk)]; \\
    1/n, & \hbox{when }\ s \in [1/(nk), \inf g^{-1}(1/n)] \\
    g(s), & \hbox{if} \ s > \inf g^{-1}(1/n),  \\
\end{array}%
\right.
$$
of continuous functions $\gamma_n$, all of them satisfying
hypothesis {\rm \bf(L)}. Obviously, $\gamma_n$ converges uniformly
to $g$ on $\R_+$. Now, for all sufficiently large $n$, Theorems
\ref{mainA} and \ref{34}  guarantee the existence of a positive
continuous function $\varphi_n(t)$ such that $\varphi_n(-\infty)
=0$, $\liminf_{t \to +\infty}\varphi_n(t) \geq \zeta_1$, and
$$
\varphi_n(t) = \frac{1}{\epsilon'
}\left\{\int_{-\infty}^te^{\lambda (t-s)}\Gamma_n(s)ds +
\int_t^{+\infty}e^{\mu (t-s)}\Gamma_n(s)ds \right\},
$$
where
$$
\Gamma_n(t):=
\int_{\R}K(s)\gamma_n(\varphi_n(t-\sqrt{\epsilon}s-h))ds.
$$
Since the shifted functions $\varphi_n(s + a)$ satisfy the same
integral equation, we can assume that $\varphi_n(0)= 0.5\zeta_1$.

Now, taking into account the inequality $|\varphi_n(t)|+
|\varphi_n'(t)|\leq \zeta_2+ \zeta_2/\epsilon', \ t \in \R,$ we
find that the set $\{\varphi_n\}$ is pre-compact in the compact
open topology of $C(\R,\R)$. Consequently we can indicate a
subsequence $\varphi_{n_j}(t)$ which converges uniformly on
compacts to some bounded element $\varphi \in C(\R,\R)$. Since
$$\lim_{j \to
+\infty}\Gamma_{n_j}(t) =
\int_{\R}K(s)g(\varphi(t-\sqrt{\epsilon}s-h))ds
$$ for every $t\in \R$, we can use the
Lebesgue's dominated convergence theorem to conclude that
$\varphi$ satisfies integral equation (\ref{iie}). Finally, notice
that $\varphi(0) = 0.5\zeta_1$ and thus $\varphi(-\infty) = 0$ (by
Theorem \ref{43}) and $\liminf_{t \to +\infty} \varphi(t) \geq
\zeta_1$ (by Theorem \ref{mainA}).

To complete the proof for Case I, we have to analyze the case when
$\max_{s \geq 0} g(s) > \max_{s \in [\zeta_1,\zeta_2]} g(s)$.
However, this cases can be reduced to the previous one if we
redefine $g(s)$ as $g(\zeta_2)$ for all $s \geq \zeta_2$, and then
observe that $\sup_{t \in \R} \varphi(s) \leq \zeta_2$ for every
solution obtained in the first part of this subsection.

\underline{\bf Case II: $c = \tilde{c}_*$}.  Consider $\epsilon_n
\uparrow 1/\tilde c_*^2$.  Then, the previous result (Case I)
assures the existence of positive functions \  $\varphi_n(t)$ \
such that $\varphi_n(-\infty) =0$, $\liminf_{t \to
+\infty}\varphi_n(t) \geq \zeta_1$, and
$$
\varphi_n(t) =
\frac{1}{\epsilon_n'}\left\{\int_{-\infty}^te^{\lambda_n
(t-s)}\Delta_n(s)ds + \int_t^{+\infty}e^{\mu_n (t-s)}\Delta_n(s)ds
\right\},
$$
where $\lambda_n < 0 < \mu_n$ satisfy $\epsilon_n z^2 -z -1 =0, \
\epsilon_n':=\epsilon_n (\mu_n - \lambda_n),$ and
$$
\Delta_n(t):= \int_{\R}K(s)g(\varphi_n(t-\sqrt{\epsilon_n}s-h))ds.
$$
The rest of proof is  exactly the same as in Case I and so is
omitted.

\vspace{0mm}

\section{Heteroclinic solutions of equation (\ref{twe})}\label{S7}
For $s \in [-\infty,0)$ and $\lambda < 0 < \mu$ satisfying
$\epsilon z^2 -z -1 =0$, set
$$
\xi(s) = \frac{\mu -\lambda}{\mu e^{-\lambda s}- \lambda e^{-\mu
s}}, \quad \mathcal{D}(s) = \min\{\int_{-h/\sqrt{\epsilon}}^{-(s
+h)/\sqrt{\epsilon}}K(u)du, \xi(-s)\}.
$$
Everywhere in this section, we assume the hypothesis
$\mathbf{(H)}$ so that all conditions of $\mathbf{(B)}$ are
satisfied. Let $\varphi(t)$ be a  semi-wavefront of equation
(\ref{twe}). Set
$$m = \liminf\limits_{t \to
+\infty}\varphi(t) \leq \limsup\limits_{t \to +\infty}\varphi(t)
=M.$$ Our next result shows that $m=M=\kappa$ if, for some $s_*
\in [-\infty,0)$, it holds
\begin{equation}\label{Dc} (1-\mathcal{D}(s_*))g'(\kappa) > -1.
\end{equation}
\begin{thm} \label{mainex} Assume {\rm \bf(H)},
$(Sg)(s)<0, s \in [\zeta_1,\zeta_2]\setminus\{s_M\},$ and
$g^2(\zeta_2) \geq \kappa$. If (\ref{Dc}) holds for some fixed
positive real number $\epsilon \leq 1/\sqrt{\tilde c_*}$, then
equation (\ref{twe}) with this $\epsilon$ has semi-wavefronts.
Moreover, each of them is in fact a wavefront.
\end{thm}

\vspace{-5mm}

\begin{pf}  Since  $\epsilon \in (0, \infty)$ and  $\tilde c_* \geq c_\#$,  Theorem \ref{mainA}
and Lemma \ref{dud} ensure that $\kappa, m, M \in
[\zeta_1,\zeta_2] $ and that $ [m,M]\subseteq g([m,M])$. The
latter inclusion and {\rm \bf(H)} imply that each of the following
three relations $\kappa\leq s_M,$ or $\kappa \leq m\leq M,$ or
$m\leq M\leq \kappa$ yields $m=M=\kappa$. Therefore, we will
consider only the case when $m < \kappa < M$ so that $g'(\kappa)
<0$. Then, by the compactness argument, we can find a solution
$y(t)$ of (\ref{twe}) such that $y(0)=\max_{s \in \R}y(s) =M$ and
$\inf_{s \in \R}y(s)\geq m$. Fix some $s_* \in [-\infty,0)$. Then
either (I) $y(t)
> \kappa$ for all $t \in [s_*,0]$ or (II) there exists some
$\hat{s} \in [s_*,0]$ such that $y(\hat{s}) = \kappa$ and $y(t) >
\kappa$ for $t \in (\hat{s},0]$.

In case (I), we have $y'(0)=0, \ y''(0) \leq 0$ and thus, in view
of Eq. (\ref{twe}),
$$
M \leq \int_{\R}K(w)g(y(-\sqrt{\epsilon}w-h))dw =
$$
$$
\int_{-h/\sqrt{\epsilon}}
^{-(s_*+h)/\sqrt{\epsilon}}K(w)g(y(-\sqrt{\epsilon}w-h))dw +
\int_{\R\setminus \mathcal{I}}K(w)g(y(-\sqrt{\epsilon}w-h))dw\leq
$$
$$
\kappa \mathcal{D}(s_*) + (1- \mathcal{D}(s_*))\max_{s \in
[m,M]}{g(s)}, \ {\rm where \ } \mathcal{I} =
[-h/\sqrt{\epsilon},-(s_*+h)/\sqrt{\epsilon}].
$$
In case (II),  considering the boundary conditions $y(\hat{s}) =
\kappa, \ y'(0)= 0$, setting
$$
G(s) = \int_{\R}K(w)g(y(s-\sqrt{\epsilon}w))dw.
$$
and then using Lemma \ref{v}, we find that \vspace{-5mm}
\begin{eqnarray*}
  M&=& y(0) = \xi (-\hat{s})\left\{\kappa +
\frac{1}{\epsilon(\mu -\lambda)}\int_{\hat{s}}^0(e^{\lambda
(\hat{s}- u)}
- e^{\mu (\hat{s}- u)})G(u-h)du\right\} \leq \\
   & & \xi (-\hat{s})\left\{\kappa + \frac{1}{\epsilon(\mu -\lambda)}\int_{\hat{s}}^0(e^{\lambda (\hat{s}- u)}
- e^{\mu (\hat{s}- u)})du\max_{x \in
[m,M]}{g(x)}\right\} =\\
  & & \hspace{-10mm} \xi (-\hat{s})\kappa +(1- \xi (-\hat{s}))\max_{s \in [m,M]}{g(s)} \leq
 \xi (-s_*)\kappa +(1-  \xi (-s_*))\max_{s \in [m,M]}{g(s)},
\end{eqnarray*}
since $\xi(-s), s \leq 0,$ is strictly increasing. Hence, we have
proved that
\begin{equation}\label{43a}
M \leq \kappa \mathcal{D}(s_*) + (1- \mathcal{D}(s_*))\max_{s \in
[m,M]}{g(s)}.
\end{equation}

\vspace{-3mm}

\noindent Analogously, there exists a solution $z(t)$ such that
$z(0)=\min_{s \in \R}z(s) =m$ and $\sup_{s \in \R}z(s)\leq M$ so
that $z'(0)=0, \ z''(0) \geq 0$. We have again that either (III)
$z(t) < \kappa$ for all $t \in [s_*,0]$ or (IV) there exists some
$\hat{s} \in [s_*,0]$ such that $z(\hat{s}) = \kappa$ and $z(t) <
\kappa$ for $t \in (\hat{s},0]$. In what follows, we are using the
condition $g^2(\zeta_2) \geq \kappa$ which implies that $g(z(t))
\geq \kappa$ once $z(t) \in [g(\zeta_2), \kappa]$. Bearing this
last remark in mind, in case (III), we obtain
$$
 m \geq \int_{\R}K(w)g(z(-\sqrt{\epsilon}w-h))dw =
$$
$$
\int_{-h/\sqrt{\epsilon}}^{-(s_*
+h)/\sqrt{\epsilon}}K(w)g(z(-\sqrt{\epsilon}w-h))dw +
\int_{\R\setminus\mathcal{ I}}K(w)g(z(-\sqrt{\epsilon}w-h))dw\geq
$$
$$
\kappa\mathcal{D}(s_*) + (1- \mathcal{D}(s_*))\min_{s \in
[m,M]}{g(s)}. $$ In case (IV), considering the boundary conditions
$z(\hat{s}) = \kappa,$ $ \ z'(0)= 0$, and using Lemma \ref{v}, we
find that \vspace{-7mm}

\begin{eqnarray*}
m&=& z(0) = \xi (-\hat{s})\left\{\kappa + \frac{1}{\epsilon(\mu
-\lambda)}\int_{\hat{s}}^0(e^{\lambda (\hat{s}- u)}
- e^{\mu (\hat{s}- u)})G(u-h)du\right\} \geq \\
   & & \xi (-\hat{s})\left\{\kappa + \frac{1}{\epsilon(\mu -\lambda)}\int_{\hat{s}}^0(e^{\lambda (\hat{s}- u)}
- e^{\mu (\hat{s}- u)})du\min_{s \in
[m,M]}{g(s)}\right\} =\\
  && \hspace{-10mm} \xi (-\hat{s})\kappa +(1- \xi (-\hat{s}))\min_{s \in [m,M]}{g(s)} \geq \xi (-s_*)\kappa +(1- \xi (-s_*))\min_{s \in [m,M]}{g(s)}.
\end{eqnarray*}
Hence, we have proved that
$$
m \geq \kappa\mathcal{D}(s_*) + (1- \mathcal{D}(s_*))\min_{s \in
[m,M]}{g(s)}.
$$
From this estimate and (\ref{43a}), we obtain that
$$
\lbrack m,M]\subseteq f([m,M])\subseteq f ^2([m,M])\subseteq
...\subseteq f ^j([m,M])\subseteq ...
$$
where $f(s) = \kappa\mathcal{D}(s_*) + (1-
\mathcal{D}(s_*)){g(s)}$ is unimodal (decreasing) if  $g$ is
unimodal (decreasing, respectively). Therefore, as $f(\kappa) =
\kappa$ and $Sf =Sg< 0$, the last chain of inclusions and the
inequality $|f'(\kappa)| \leq 1$ is sufficient to obtain
$m=M=\kappa$, see Proposition \ref{SW}. \hfill {$\blacksquare$}
\end{pf}

\vspace{-5mm}

\begin{rem}Theorem \ref{mainex2} follows from Theorem \ref{mainex} if we
take $s_* = -h$ and observe that $0 < e^{-h} \leq \xi(h)  < 1$. Note
that $e^{-h} \leq \xi(h)$ amounts to the inequality $ \mu (1 -
e^{-h(\lambda+1)}) \geq \lambda (1 - e^{-h(\mu+1)}), $ which holds
true since the left hand side is positive and the right hand side is
negative.
\end{rem}
\begin{rem}
For fixed $\epsilon,h$, and for $s < 0$, consider the following
equation
\begin{equation}\label{dela}
\int_{-h/\sqrt{\epsilon}}^{-(s +h)/\sqrt{\epsilon}}K(u)du =
\xi(-s).
\end{equation}
It is clear that the left hand side of (\ref{dela}) is decreasing
in $s \in (-\infty,0)$ from
$\int_{-h/\sqrt{\epsilon}}^{+\infty}K(u)du \geq 0$ to $0$  while
the right hand side is strictly increasing from $0$ to $1$. Thus
(\ref{dela})  has a unique  solution $s' \in [-\infty,0)$ which
coincides with the optimal value of $s_*$ in (\ref{Dc}).
\end{rem}

\vspace{-4mm}

\section{Proof of Theorem \ref{main} }
\label{rrr}  For the convenience of the reader, the proof will be
divided in several steps. Note that the assumptions of Theorem
\ref{main} imply that supp $K\cap (-h/\sqrt{\epsilon}, \eta)\not=
\emptyset$.

\vspace{-2mm}

{\it Claim I:  $y(t):=\varphi(t) - \kappa > 0$ is not
superexponentially small as $t \to +\infty$.}\\
Let $\varphi: \R \to (0,+\infty)$ be a non-constant solution of
(\ref{twe}) satisfying $\varphi(+\infty)=\kappa$. First we prove
that $\varphi$ cannot be eventually constant. Indeed, if
$\varphi(t) = \kappa$ for all $t \geq -h $ and $\varphi(t)$ is not
constant in some left neighborhood of $t =-h$ then  we obtain from
(\ref{twe}) that
\begin{equation}\label{pe}
\int_{-\eta}^{\eta}K(s)q(t-\sqrt{\epsilon}s)ds \equiv \kappa,
\quad t \in [-h, \sqrt{\epsilon}\eta],
\end{equation}
where $q(t)=g(\varphi(t-h))$. Set $K_1(u) = K(-u/\sqrt{\epsilon}
+\eta)/\sqrt{\epsilon},$ $p(u) = q(-u) - \kappa$, $ \ x =
\sqrt{\epsilon}\eta-t, \ t \in [-h, \sqrt{\epsilon}\eta]$. Then
(\ref{pe}) can be written as  a scalar Volterra convolution
equation on a finite interval
$$
\int_{0}^{x}K_1(x-s)p(s)ds
=\int_{t/\sqrt{\epsilon}}^{\eta}K(s)(q(t-\sqrt{\epsilon}s)-\kappa)ds
\equiv 0, \ x \in [0, \sqrt{\epsilon}\eta+h].
$$
In consequence, since supp $K\cap (-h/\sqrt{\epsilon}, \eta)\not=
\emptyset$, a result of Titchmarsh (see \cite[Theorem 152]{tit})
implies that
$$p(s)= g(\varphi(-s-h)) - \kappa = 0, \quad s \in [0,
\sqrt{\epsilon}\eta+h].$$ Thus $\varphi(t)= \kappa$ for all $t \in
[-2h-\sqrt{\epsilon}\eta,-h]$, a contradiction.

Now, when $\varphi$ is not oscillating around the positive
equilibrium, we can see that $y(t) = \varphi(t)-\kappa$ is either
decreasing and strictly positive or increasing and strictly
negative, for all sufficiently large $t$.  Indeed, if $\varphi(t)
\geq \kappa, \ t \geq -h - \sqrt{\epsilon} \eta,$ has a local
maximum at $t=b > 0$ then $\varphi(b) > \kappa, \ \varphi'(b) =0,
\ \varphi''(b) \leq 0$. In consequence, since
$\varphi(+\infty)=\kappa$ and $g'(\kappa) <0$, we get
$$
\kappa < \varphi(b) \leq
\int_{-\eta}^{\eta}K(s)g(\varphi(b-\sqrt{\epsilon}s-h))ds \leq
\int_{-\eta}^{\eta}K(s)g(\kappa)ds = \kappa,
$$
a contradiction. The same argument works when $\varphi(t) \leq
\kappa$ for all large $t$.

Next, observe that $y(t)$ satisfies $ \epsilon y''(t) - y'(t)=
y(t)+k(t)y(t-\sqrt{\epsilon}\eta-h), $ where, in view of the
monotonicity of $y$, it holds that $-2g'(\kappa) \geq k(t):=$
$$
=
-\int_{-\eta}^{\eta}K(s)\frac{g(\varphi(t-\sqrt{\epsilon}s-h))-g(\kappa)}{\varphi(t-\sqrt{\epsilon}s-h)-\kappa}\cdot
\frac{\varphi(t-\sqrt{\epsilon}s-h)-\kappa}{\varphi(t-\sqrt{\epsilon}\eta-h)-\kappa}ds
\geq 0
$$
for all sufficiently large $t$. We can use now Lemma 3.1.1 from
\cite{hl2} to conclude that $y(t)> 0$ cannot converge
superexponentially to $0$.\\
{\it Claim II: $y(t)> 0$ cannot hold when Eq. (\ref{cchh1}) does
not have roots in $(-\infty,0)$}. Observe that
$y(t)=\varphi(t)-\kappa,$ $ y(+\infty) = 0,$ verifies
$$\epsilon y''(t) - y'(t)- y(t)+
\int_{-\eta}^{\eta}K(s)g_1(y(t-\sqrt{\epsilon}s-h))ds=0,\ t \in
\R,\
$$
where $g_1(s):= g(s+ \kappa)- \kappa, \ g_1(0)=0, g_1'(0)=
g'(\kappa)<0$.  In virtue of Claim I  and Lemma \ref{41aa}, we can
find a real number $d
>1$ and a sequence $t_n \to +\infty$ such that $y(t_n) = \max_{s
\geq t_n}y(s)$ and
$$
\max_{s \in [t_n- 3h- 3\eta\sqrt{\epsilon},t_n]}y(s) \leq d
y(t_n).
$$
Additionally, we can find a sequence $\{s_n\},\ \lim (s_n-t_n) =
+\infty$ such that $|y'(s_n)| \leq y(t_n)$. Now, $w_n(t) =
y(t+t_n)/y(t_n), \ t \in \R $ satisfies
$$
\epsilon w''(t) - w'(t)-w(t)+
\int_{-\eta}^{\eta}K(s)p_n(t-\sqrt{\epsilon}s-h)w(t-\sqrt{\epsilon}s-h)ds=0,
$$
where $p_n(t) = g_1(y(t+t_n))/y(t+t_n)$. It is clear that $\lim
p_n(t) =g'(\kappa)$ for every $t \in \R,$ and that $0 < w_n(t)
\leq d$ for all $t \geq -3(\eta\sqrt{\epsilon} +h)$.

To estimate $|w'_n(t)|$,  let $\mathscr{W}_n(t):=
\int_{-\eta}^{\eta}K(s)p_n(t-\sqrt{\epsilon}s-h)w_n(t-\sqrt{\epsilon}s-h)ds$.
Since $z_n(t) = w'_n(t)$ satisfies $z_n(s_n-t_n) = y'(s_n)/y(t_n)
\in [-1,0]$ and
$$
\epsilon z_n'(t) - z_n(t)-w_n(t)+\mathscr{W}_n(t)=0, \quad t \in
\R ,
$$
\vspace{-3mm} we obtain that \begin{equation}\label{whi} w'_n(t) =
e^{(t+t_n-s_n)/\epsilon}z_n(s_n-t_n) +
\frac{1}{\epsilon}\int_{s_n-t_n}^te^{(t-s)/\epsilon}(w_n(s)-
\mathscr{W}_n(s))ds.
\end{equation}
\vspace{-2mm} Furthermore, for each fixed $t\geq
-2\eta\sqrt{\epsilon}-2h$ and sufficiently  large $n$, we have
$$|w'_n(t)|\leq 1+
\frac{1}{\epsilon}\int^{s_n-t_n}_te^{(t-s)/\epsilon}(w_n(s)+\sup_{s\not=0}\frac{|g_1(s)|}{|s|}
\int_{-\eta}^{\eta}K(u)w_n(s-\sqrt{\epsilon}u-h)du)ds\leq
$$
$$\leq 1+ (\sup_{s\not=0}\frac{|g_1(s)|}{|s|}+d)
\frac{1}{\epsilon}\int^{s_n-t_n}_te^{(t-s)/\epsilon}ds\leq 1+d+
\sup_{s\not=0}\frac{|g_1(s)|}{|s|}.
$$
Hence, there is a subsequence $\{w_{n_j}(t)\}$ which converges on
$[-2\eta\sqrt{\epsilon}-2h, +\infty)$, in the compact-open
topology, to a non-negative decreasing function $w_*(t), $ $ \
w_*(0) = 1,$ such that $w_*(t) \leq d$ for all $t \geq
-2\eta\sqrt{\epsilon}-2h$. By the Lebesgue bounded convergence
theorem, we find, for all $t \in [-\eta\sqrt{\epsilon}-h,
+\infty)$, that $$\mathscr{W}_{n_j}(t) \to
g'(\kappa)\int_{-\eta}^{\eta}K(s)w_*(t-\sqrt{\epsilon}s-h)ds.
$$
In consequence, integrating (\ref{whi}) between $0$ and $t$ and
then taking the limit as $n_j \to \infty$ in the obtained
expression, we establish that $w_*(t)$ satisfies
\begin{equation}\label{eqwo}
\epsilon w''(t) - w'(t)-w(t)+
g'(\kappa)\int_{-\eta}^{\eta}K(s)w(t-\sqrt{\epsilon}s-h)ds=0
\end{equation}
for all $t \geq -\eta\sqrt{\epsilon}-h$. We claim that $w_*(t)$ is
positive for $t \geq -\eta\sqrt{\epsilon}-h$. Indeed, if $w_*(t')
=0$ for some $t'$ then $t'>0$ since $w_*(0) =1$ and $w_*(t)$ is
decreasing. Next, if $t'$ is the leftmost positive point where
$w_*(t') =0$, then (\ref{eqwo}) implies
$$
\int^{\eta}_{\max\{-\eta,-h/\sqrt{\epsilon}\}}K(s)w(t'-\sqrt{\epsilon}s-h)ds=0.
$$
However, this contradicts to the  following  two facts: (i) due to
the definition of $t'$, it holds $w(t'-\sqrt{\epsilon}s-h) > 0$
for all $s \in (\max\{-\eta,-h/\sqrt{\epsilon}\}, \eta)$; (ii)
$K(s) \geq 0$ and supp $K\cap (-h/\sqrt{\epsilon}, \eta)\not=
\emptyset$. Hence, $w_*(t)
>0$ and we can use Lemma 3.1.1 from \cite{hl2} to conclude
that $w_*(t)> 0$ is not a small solution. Then Lemma \ref{41lin}
implies that there exists $\gamma <0$ such that $w_*(t) = v(t) +
O(\exp(\gamma t)), \ t \to +\infty,$ where $v$ is a {\it non
empty} finite sum of eigensolutions of (\ref{eqwo}) associated to
the eigenvalues $\lambda_j \in F= \{\gamma < \Re\lambda_j \leq
0\}$. Now, since the set $F$ does not contain any real eigenvalue
by our assumption, we conclude that $w_*(t)$ should be oscillating
on $\mathbb{R}_+$ (e.g. see \cite[Lemma 2.3]{hl1}),  a
contradiction. \hfill {$\blacksquare$}
\begin{rem}To establish the non-monotonicity of wavefronts in
\cite{FT}, the hyperbolicity of equation (\ref{cchh1}) and
$C^2$-smoothness of $g$ at $\kappa$ were assumed. However, as we
have shown, the first condition can be removed and it suffices to
assume that $g$ is a continuous function which is differentiable
at $\kappa$.
\end{rem}
\begin{rem} \label{re} For equation (\ref{17}), Liang and Wu
found numerically that the wavefronts may exhibit unsteady
multihumps. As it is observed in \cite{gouss,LiWu} for these
cases, the first hump (its shape, size and location) remains
stable on the front of the waves, but the second hump expands in
width to the positive direction as the number of iteration is
increasing. However the multihump waves of \cite{gouss,LiWu} may
appear due to the numerical instability of the algorithms. Indeed,
assuming {\bf (H)} and reasoning as in Section \ref{S7}, we find
that, for a fixed $\alpha > \kappa$, neither wavefront $\phi(t)$
can satisfy $\phi(t) \geq \alpha$ during 'sufficiently large'
period of time $J$ (the maximal admissible length of $J$ depends
on $\alpha$: $|J|= 2q_*(\alpha)
>0$). Indeed, supposing that $M = y(0)= \max_{s \in J} y(s)$ and that $q_*$ is sufficiently large, we
 get a contradiction:
$$ M
\leq \int_{\R}K(w)g(y(-\sqrt{\epsilon}w-h))dw =
$$
\vspace{-5mm}
$$
\int_{-(q_*+h)/\sqrt{\epsilon}}^{(q_*
-h)/\sqrt{\epsilon}}K(w)g(y(-\sqrt{\epsilon}w-h))dw +
\int_{\R\setminus J}K(w)g(y(-\sqrt{\epsilon}w-h))dw\leq
$$
$$
g(\alpha) \mathcal{D}_1(q_*) + (1- \mathcal{D}_1(q_*))\max_{x \geq 0
}{g(x)} < \kappa < M,
$$
since
$$
\lim_{q_* \to +\infty} \mathcal{D}_1(q_*):= \lim_{q_* \to +\infty}
\int_{-(q_*+h)/\sqrt{\epsilon}}^{(q_* -h)/\sqrt{\epsilon}}K(w)dw =1.
$$
\end{rem}

\vspace{-5mm}
\section{Appendix}

\vspace{-3mm}

Consider $ \psi(z,\epsilon) =\epsilon z^2-z-q+p\exp(-zh)\int_ {\R}
K(s)\exp(-\sqrt{\epsilon} zs)ds, $ where  $p> q$ and $K(s)$
satisfies condition (\ref{kkk}).
\begin{lem} \label{L23} Assume that $p > q >0$.
Then there exist extended positive real numbers $\epsilon_0 \leq
\epsilon_1, \epsilon_i = \epsilon_i(h,p,q),$ such that, for every
$\epsilon \in (0, \epsilon_0)\cup (\epsilon_1, \infty)$, equation
$\psi(\lambda,\epsilon) =0$ has exactly two real roots
$\lambda_1(\epsilon)<\lambda_2(\epsilon)$. Furthermore,
$\lambda_1(\epsilon), \lambda_2(\epsilon)$ are positive if
$\epsilon < \epsilon_0$ and are negative if $\epsilon >
\epsilon_1$. If $\epsilon \in (\epsilon_0, \epsilon_1)$, then
$\psi(z,\epsilon)> 0$ for all $z \in \R$. Next, $\epsilon_0 =
\epsilon_1$  if and only if $\epsilon_0=\epsilon_1 = \infty$.
Furthermore, $\epsilon_1 =\infty$ if $\int_{\R}sK(s)ds \geq 0$ and
$\epsilon_1 $ is finite if the equation
\begin{equation}\label{nega} z^2 -q + p\int_{\R}\exp(-zs)K(s)ds=0
\end{equation} has two
negative roots. Finally, if $\int_{\R}xK(x)dx \leq 0$ then
$\epsilon_0$ is finite  and
\begin{equation}\label{am}
c_*: = 1/\sqrt{\epsilon_0} > |\int_{\R}sK(s)ds|/(h+ 1/p).
\end{equation}
\end{lem}

\vspace{-10mm}
\begin{pf}  Observe that $\psi''_z(z,\epsilon) > 0, \ z \in \R$, so that
$\psi(z,\epsilon)$  is strictly concave with respect to $z$.  This
guaranties the existence of at most two real roots. Next, since
$\psi(z,0)$ has a unique real (positive) root $z_0$, where
$\psi'_z(z_0,0)<0$, we find that $\psi(z,\varepsilon)$ possesses
exactly two positive roots for all small $\epsilon >0$.

After introducing a new variable $w = \sqrt{\epsilon}z$, we find
that equation $\psi(z,\epsilon) =0$ takes the following form
\begin{equation}\label{tr}
(q+\frac{w}{\sqrt{\epsilon}}-w^2)\exp({\frac{wh}{\sqrt{\epsilon}}})=
p\int_{\R}\exp(-ws)K(s)ds \quad(:= G(w)).
\end{equation}
As we have seen,  equation (\ref{tr}) may have at most two real
roots and, for small $\epsilon >0$, it possesses two positive
roots $w_1(\epsilon) < w_2(\epsilon)$. Furthermore, we have that
$G(0)=p, \  G''(w) > 0$. An easy analysis of (\ref{tr}) shows that
positive $w_1(\epsilon)< w_2(\epsilon)$ exist and depend
continuously on $\epsilon$ from the maximal interval $(0,
\epsilon_0)$, where $\epsilon_0$, when finite,  is determined by
the relation $w_1(\epsilon_0) = w_2(\epsilon_0)$. To prove that
equation (\ref{tr}) does not have any real positive root for
$\epsilon > \epsilon_0$, it suffices to note that $G(w)$ does not
depend on $\epsilon$ while the left hand side of (\ref{tr})
decreases with respect to $\epsilon$ at every positive point $w$
where $q+w/\sqrt{\epsilon}-w^2 >0$.

Similarly, for $\epsilon > \epsilon_0$,  the left hand side of
(\ref{tr}) increases to $q - w^2$ with respect to $\epsilon$ at
every $w<0$ where $q+w/\sqrt{\epsilon}-w^2 >0$. In consequence,
$\epsilon_1 $ is finite if and only if equation (\ref{nega}) has
two simple negative roots. It is evident that this may happen only
if $G'(0) = -\int_{\R}sK(s)ds > 0$ and that in this case
$\epsilon_1 > \epsilon_0$.

Clearly,  $\psi'_z(0,\epsilon_0) < 0$. For $\int_{\R}xK(x)dx \leq
0$, the latter inequality amounts to (\ref{am}). It is easy to see
that the equality $\epsilon_0 = +\infty$ actually can happen when
$\int_{\R}xK(x)dx
> 0$. \hfill {$\blacksquare$}
\end{pf}

\vspace{-5mm}

Next propositions are crucial in the proof of Theorem \ref{mainex}.

\begin{lem} \label{41aa} Let $x: \R_+ \to (0, +\infty)$ satisfy
$x(+\infty)=0$. Given an integer $d >1$ and a real $\rho > 0$, we
define $\alpha = (\ln d)/\rho > 0$. Then either (a) $x(t) =
O(e^{-\alpha t})$ at $t = +\infty$, or (b) there exists a sequence
 $t_j \to +\infty$ such that $x(t_j) = \max_{s \geq t_j}x(s)$ and
 $\max_{s \in [t_j-\rho,t_j]}x(s) \leq d x(t_j)$.
\end{lem}

\vspace{-5mm}
\begin{pf} Set
$ T= \left\{ t: x(t) = \max_{s \geq t}x(s)  \ {\rm and} \ \max_{s
\in [t-\rho,t]}x(s) \leq d x(t) \right\} . $ \\ Then either (I) $T
\not= \emptyset$ and $\sup T = +\infty$ and therefore the
conclusion (b) of the lemma holds, or (II) $T$ is a bounded set
(without restricting the generality, we may assume that $T=
\emptyset$). Let us analyze more closely the second case
(supposing that $T= \emptyset$). Take an arbitrary $t > \rho$ and
let $\hat{t} \geq t$ be defined as leftmost point where
$x(\hat{t})= \max_{s \geq t}x(s)$ . Since $t \not\in T$, we have
that $\hat{t}-t \leq \rho$. Let $t_1$ be defined by $x(t_1) =
\max_{s \in [\hat{t}-\rho, \hat{t}]}x(s)$, our assumption about
$T$ implies that $\hat{t}-\rho \leq t_1 < t \leq \hat{t}$ and that
$x(t_1)
> dx(\hat{t}) \geq dx(t)$. Additionally, $x(t_1)= \max_{s \geq
t_1}x(s)$. Next, we define $t_2$ as leftmost point satisfying
$x(t_2) = \max_{s \in [t_1-\rho, t_1]}x(s)$. Notice that $0 <
t_1-t_2 \leq \rho$ and $x(t_2) > dx(t_1)$. Proceeding in this way,
we construct a decreasing sequence $t_j$ such that $x(t_{j+1}) >
dx(t_{j})$ for every $j$. We claim that there exist an integer $m$
such that $t_m \leq \rho$. Indeed, otherwise $t_j> \rho$ for all $j
\in \N$ that implies the existence of $\lim t_j = t_*$ and $\lim
x(t_j) = x(t_*)$. However, this is not possible since $x(t_j)
> d^jx(t)>0$. Hence, $t_m \in [0,\rho]$ for some integer $m$.
Notice that $t-t_m\leq \hat{t}-t_m \leq m\rho$ implying that $m \geq
(t-t_m)/\rho \geq (t-\rho)/\rho$ and that
$$x(t) < d^{-m}x(t_m) \leq d^{-m}\max_{s
\in [0,\rho]}x(s) \leq de^{-\alpha t}\max_{s \in [0,\rho]}x(s).
\qquad {\blacksquare}$$
\end{pf}
The proof of the next lemma follows  that of Proposition 7.1 from
\cite{FA}. When $K(s)=\delta(s)$ is a Dirac delta function, the
obtained asymptotic estimates for $x$ are uniform in $\epsilon$,
see \cite[Lemma 4.1]{AT}.

\vspace{0mm}
\begin{lem} \label{41lin}  Let
$x\in C^2(\R,\R)$ verify the  equation
\begin{equation}\label{twel}
x''(t) +\alpha x'(t)+\beta x(t)+ p\int_{\R}K(s)x(t+qs+h)ds = f(t), \
t \geq 0,
\end{equation}
where $K$ satisfies (\ref{kkk}), $\alpha,\beta, p,q,h \in \R$ and
$f(t) = O(\exp(-bt)), \ t \to +\infty$ for some $b> 0$.  Suppose
further that $|x(t)| \leq c \exp({\gamma t}), \ t \leq 0,$ for
some $\gamma \leq 0$, and that $ \sup_{t \geq 0} |x(t)|$ is
finite. Then, given $\sigma \in (0,b)$, it holds that
$$
x(t) = w(t) + \exp(-(b-\sigma)t)o(1), \ t \to +\infty,
$$
where $w(t)$ is a finite sum of  eigensolutions of (\ref{twel})
associated to the eigenvalues $\lambda_j \in \{-(b-\sigma) <
\Re\lambda_i \leq 0\}$.
\end{lem}
\vspace{-5mm}
\begin{pf} Remark that the conditions of Lemma \ref{41lin} imply that
$\sup_{t \geq 0} |x''(t)|$ is finite and that $|x'(t)| = O(1)$ at
$t =+\infty$ (if $\alpha\not=0$) or $|x'(t)| = O(t)$ (if $\alpha
=0$). The proof of this observation is based on deriving
estimations similar to (\ref{nv}) and is omitted here. Applying
the Laplace transform $\mathcal{L}$ to (\ref{twel}), we obtain
that $ \chi(z)\tilde{x}(z)= \tilde{f}(z) + r(z), $ where
$\tilde{x} = \mathcal{L}x, \ \tilde{f} = \mathcal{L}f$ and
$$
r(z) = x'(0)+zx(0)+\alpha x(0)+pe^{zh}\int_{\R}K(s)
e^{zqs}ds\int^{h+ps}_0 e^{-zu}x(u)du,
$$
$$
\chi(z) =z^2+\alpha z + \beta+pe^{zh}\int_ {\R} K(s)e^{qzs}ds.
$$
Since  $x$ is bounded on $\R_+$, we conclude that $\tilde{x}$ is
analytic in $\Re z
> 0$. Moreover, from the growth restrictions 
on $x, f, K$ we obtain that $r$ is an entire function and
$\tilde{f}$ is holomorphic in $\Re z > -b$. Therefore
$H(z)=(\tilde{f}(z) + r(z))/\chi(z)$ is meromorphic in $\Re z
> -b$. Observe also that $H(z) = O(z^{-1}),\ z \to \infty,$
for each fixed strip $\Pi(s_1,s_2)=\{s_1 \leq \Re z \leq s_2\},
s_1> -b$. Now, let $ \sigma > 0$ be such that the vertical strip
$-b < \Re z < -b + 2\sigma$ does not contain any zero of
$\chi(z)$. By the inversion formula, for some sufficiently small
$\delta >0$, we obtain that
$$
x(t) = \frac{1}{2\pi i}\int_{\delta -i\infty}^{\delta
+i\infty}e^{zt}\tilde{x}(z)dz = \frac{1}{2\pi i}\int_{\delta
-i\infty}^{\delta +i\infty}e^{zt}H(z)dz =w(t) + u(t),  \ t
>0,
$$
$$
\hspace{-11mm}{\rm where} \ w(t) = \sum_{-b + \sigma < \Re
\lambda_j \leq 0} {\rm Res}_{z=\lambda_j}
\frac{e^{zt}(\tilde{f}(z) + r(z))}{\chi(z)} = \sum _{-b +\sigma <
\Re \lambda_j \leq 0} e^{\lambda_jt}P_j(t) ,
$$
$$
u(t) = \frac{1}{2\pi i}\int_{-b +\sigma -i\infty}^{-b +\sigma
+i\infty}e^{zt}H(z)dz.
$$
Now, observe that on any vertical line in $\Re z
> -b$ which does not pass through
the poles of $\chi(z, \epsilon)$ and $0 \in \mathbb{C}$, we have
$$H(z) = a(z) + \frac{x(0)}{z}, \ {\rm where \ } a(z) = O(z^{-2}), \
z \to \infty.$$ Therefore, for $a_1(s) = a(-b+\sigma + is)$, we
obtain
$$
u(t) =\frac{e^{(-b + \sigma)t}}{2\pi i}\left\{
\int_{-\infty}^{+\infty}e^{ist}a_1(s)ds \right\} + \frac{x(0)}{2\pi
i}\int_{-b +\sigma -i\infty}^{-b +\sigma +i\infty}\frac{e^{zt}}{z}dz
, \ t > 0.$$ Next, since $a_1 \in L_1(\R)$, we have, by the
Riemann-Lebesgue lemma, that
$$
\lim_{t \to +\infty}\int_{\R}e^{ist}a_1(s)ds =0.
$$
For $t >0$, a direct computation shows that $\int_{-b +\sigma
-i\infty}^{-b +\sigma +i\infty}z^{-1}e^{zt}dz =0.$ Thus we get
$u(t) = e^{-(b - \sigma)t} o(1)$, and the proof is completed.
\hfill {$\blacksquare$}
\end{pf}
\begin{lem} \label{v} If $x$ verifies (\ref{twe}) and the
conditions $x(a) = x_0,  x'(b) =0$, then
$$
x(b) = \xi(b-a)\left\{x_0 + \frac{1}{\epsilon(\mu
-\lambda)}\int_{a}^b(e^{\lambda(a- u)} - e^{\mu(a-
u)})(Gx)(u-h)du\right\}.
$$
\end{lem}
\vspace{-5mm}
\begin{pf} It suffices to consider the variation of constants
formula for (\ref{twe}):
$$
x(t) = Ae^{\lambda t} + Be^{\mu t} + \frac{1}{\epsilon (\mu -
\lambda)}\left\{ \int_a^te^{\lambda (t-s)}g(s)ds + \int_t^be^{\mu
(t-s)}g(s)ds \right\},
$$
where $g(s):= (Gx)(s-h)$. \hfill {$\blacksquare$}
\end{pf}
\vspace{-3mm} The following proposition is due to Singer
\cite{LMT}:
\begin{prop}\label{SW}
 Assume that $f: [\zeta_*, \zeta^*] \to [\zeta_*, \zeta^*],
 \  f \in C^3[a, b]$, is either strictly decreasing function or it has only one
critical point $x_M$ (maximum) in $[\zeta_*, \zeta^*]$. If the
unique fixed point $\kappa \in [\zeta_*, \zeta^*]$ of $f$ is locally
asymptotically stable and the Schwarzian derivative satisfies $(S
f)(s) < 0$ for all $s \not=s_M$ then $\kappa$ is globally
asymptotically stable.
\end{prop}
The condition of the negativity of $Sg$ (which requires $C^3-$
smoothness of $g$) can be  weakened with the use of a generalized
Yorke condition introduced in \cite{ltt} and analyzed in
\cite{bio} from the biological point of view.

\vspace{-3mm}

\section*{Acknowledgments}

\vspace{-3mm}

 This work was supported by CONICYT (Chile) through
PBCT program ACT-05 and FONDECYT (Chile) projects 1030992 (P.A.),
1071053 (S.T.), and by the University of Talca, program
``Reticulados y Ecuaciones". E. Trofimchuk was supported by DFFD
(Ukraine) project $\Phi$25.1/021. We are grateful to  Gabriel
Valenzuela for suggesting the present version of Lemma \ref{L23}. We
would like to thank the anonymous referee for the constructive
criticism and useful comments which helped us to improve
substantially the paper.

\end{document}